\documentclass{article}

\usepackage{amssymb,latexsym,amsfonts}

\title{Degenerate Gauss hypergeometric functions}

\author{Raimundas Vid\=unas\footnote{Supported by the 21 Century
COE Programme "Development of Dynamic Mathematics with High Functionality"
of the Ministry of Education, Culture, Sports, Science and Technology of
Japan.}\\ \em Kyushu University}

\newtheorem{theorem}{Theorem}[section]
\newtheorem{lemma}[theorem]{Lemma}
\newtheorem{conclude}[theorem]{Corollary}

\usepackage{latexsym}

\newcommand{\hpg}[5]{{}_{#1}\mbox{\rm F}_{\!#2}\!
  \left(\left.{#3 \atop #4}\right| #5 \right) }

\newcommand{\hpgo}[2]{{}_{#1}\mbox{\rm F}_{\!#2}}

\newcommand{\hppg}[3]{%{}_2\mbox{\bf F}_{\!1}\!
  \mbox{\bf F}\!\left(\left.{#1 \atop #2}\right| #3 \right)}

\newcommand{\proof}{{\bf Proof. }}
\newcommand{\qed}{\hfill \mbox{$\Box$}\\}
\newcommand{\equal}{\!\!\!=\!\!\!}
\newcommand{\CC}{{\Bbb C}}

\newcommand{\PP}{{\Bbb P}}
\newcommand{\GG}{{\Bbb G}}
\newcommand{\ZZ}{{\Bbb Z}}
\newcommand{\ZZnp}{\ZZ_{\le 0}}
\newcommand{\ZZnn}{\ZZ_{\ge 0}}

\begin{document}

\maketitle

\begin{abstract} This paper studies terminating and ill-defined Gauss hypergeometric
functions. For their hypergeometric equations, the set of 24 Kummer's
solutions degenerates. We describe those solutions and relations between
them.
\end{abstract}

\section{Introduction}
\label{logarithms}

Throughout the paper, let $\ZZnp$ and $\ZZnn$ denote the sets %$\{0,-1,-2,\ldots\}$
of non-positive and non-negative integers, respectively. The {\em Gauss
hypergeometric function} is defined by the series
\begin{equation} \label{hpggauss}
\hpg{2}{1}{\!a,\,b\,}{c}{z}=1+\frac{a\;b}{c\!\cdot 1!}\,z+
\frac{a\,(a\!+\!1)\,b\,(b\!+\!1)}{c\,(c\!+\!1)\cdot 2!}\,z^2 +\ldots.
\end{equation}
The hypergeometric series terminate %and represent a polynomial function
when $a\in\ZZnp$ or \mbox{$b\in\ZZnp$}; %are non-positive integers.
the function is undefined when $c\in\ZZnp$. If $n,N$ are non-negative integers %\in\ZZnn$
and $n\le N$, then we interpret (\ref{hpggauss}) as a terminating
hypergeometric series. % as well, with at most $n$ terms.
Gauss hypergeometric function can be analytically continued onto
$\CC\setminus(0,\infty)$. It satisfies the {\em hypergeometric differential
equation}
\begin{equation} \label{hpgde}
z\,(1-z)\,\frac{d^2y(z)}{dz^2}+
\big(c-(a\!+\!b\!+\!1)\,z\big)\,\frac{dy(z)}{dz}-a\,b\,y(z)=0.
\end{equation}
Throughout the paper, we denote this equation by $E(a,b,c)$. Since Kummer we
know that in general there are 24 hypergeometric series which express
solutions of $E(a,b,c)$. %In general, the same Gauss hypergeometric function
%can be represented by 4 distinct hypergeometric series around $z=0$.

The subject of this paper is Gauss hypergeometric functions when some or all
of the numbers $a$, $b$, $c$, $c-a$, $c-b$, $a-b$, $c-a-b$ are integers. We
refer to these functions and to corresponding hypergeometric equations as
{\em degenerate}. (Note that in \cite{bateman} and \cite{abrostegun}, a
hypergeometric equation is called degenerate only when $a$, $b$, $c-a$ or
$c-b$ are integers.) In a degenerate case, some of the 24 Kummer's
hypergeometric solutions are terminating or undefined, relations between
them degenerate. The monodromy group of $E(a,b,c)$ deserves a separate
description in most of these cases. In particular, there are terminating
hypergeometric solutions if and only if the monodromy group is reducible.

Bases of solutions for degenerate hypergeometric equations are presented in
\cite[Sections 2.2-2.3]{bateman}. Some of those solutions are reproduced in
\cite[Section 15.5]{abrostegun}. Logarithmic solutions are derived in
several texts, for example in \cite[pg. 82-84]{specfaar}. Generators for the
monodromy group in all degenerate cases are given in \cite[Section
4.3]{gausspainleve}. Transformations of terminating Gauss hypergeometric
series are widely used (see \cite[Section 0.6]{koekswart} for example), but
they are seldom fully exhibited.
% See also \cite[Section 15.5]{abrostegun}. \cite[Section 2.2]{specfaar}
Here we present hypergeometric solutions of degenerate hypergeometric
equations and relations between them more explicitly and systematically.
Our Theorem \ref{dendcases} and Corollary \ref{abcases} %and Table \ref{figtab}
give convenient characterization of various degenerate cases.

\section{General observations}

Hypergeometric equation $E(a,b,c)$ is a Fuchsian equation with three regular
singular points $z=0$, $z=1$ and $z=\infty$. The local exponents are:
\[
\mbox{$0$, $1-c$ at $z=0$;} \qquad \mbox{$0$, $c-a-b$ at $z=1$;}\qquad
\mbox{and $a$, $b$ at $z=\infty$.}
\] %differences at these points are $1-c$, $c-a-b$, $a-b$ respectively.
Permutations of the three singular points %$z=0,z=1,z=\infty$
and of their local exponents are realized by transformations $y(z)\mapsto
z^{-\alpha}(1-z)^{-\beta}\,y(\varphi(z))$ of $E(a,b,c)$ to other
hypergeometric equations, where
\begin{equation} \label{p1automor}
\varphi(z)\in\left\{ z,\,\frac{z}{z-1},\,1-z,\,1-\frac1z,\,\frac1z,\,
\frac1{1-z}\right\}
\end{equation}
and $\alpha,\beta$ are suitably chosen from the set $\{0,a,b,1-c,c-a-b\}$.
We refer to these transformations of hypergeometric equations as the {\em
fractional-linear transformations}.

In general, $E(a,b,c)$ can be transformed to 23 other hypergeometric
equations by the fractional-linear transformations. Consequently, there are
24 hypergeometric series representing solutions of $E(a,b,c)$; they are
commonly referred to as the {\em $24$ Kummer's solutions}. These solutions
and general relations between them are fully presented in \cite[Section
2.9]{bateman}. In general, the 24 hypergeometric series represent 6
different Gauss hypergeometric functions, since Euler's and Pfaff's formulas
\cite[Theorem 2.2.5]{specfaar} identify four Kummer's series with each
other. We refer to those identities as the {\em Euler-Pfaff
transformations}.

If we consider permutation of the upper parameters $a$, $b$ as a non-trivial
transformation, we have a group of 48 elements acting on hypergeometric
equations. This group acts on the parameters of $E(a,b,c)$ as follows.
\begin{lemma} \label{frlints}
\begin{enumerate}
\item Fractional-linear transformations can permute the three numbers
\mbox{$1-c$}, $c-a-b$, $b-a$ and change their signs in any way.%
\item Fractional-linear transformations can permute the four numbers%
\[ %begin{equation} \label{constperm}
\textstyle -\frac12+a,\;\frac12-b,\;-\frac12+c-a,\;\frac12+b-c,
\] %end{equation}
in any way, and can change their signs simultaneously.
\end{enumerate}
\end{lemma}
\proof We define $e_0=1-c$, $e_1=c-a-b$, $e_\infty=b-a$. Note that $e_0$,
$e_1$, $e_\infty$ are the local exponent differences of $E(a,b,c)$ at the
singular points. The first statement is clear once one accepts
characterization of fractional-linear transformations as the transformations
of hypergeometric equations which permute the three singular points and
their local exponents. (Notice that interchanging local exponents at a
singular point changes the sign of the local exponent difference.)

For the second statement, notice that the four listed numbers are equal to
\[
\frac{-e_0-e_1-e_\infty}2,\qquad \frac{e_0+e_1-e_\infty}2,\qquad
\frac{-e_0+e_1+e_\infty}2,\qquad \frac{e_0-e_1+e_\infty}2.
\]
By the first statement, we can permute $e_i$'s and change their signs. If we
change the signs of even number of $e_i$'s, the four numbers get permuted;
otherwise they get multiplied by $-1$ and permuted.
%Quite evidently, this claim follows from the first statement.

Alternatively, one can show both statements by checking the list of 24
related hypergeometric equations:\vspace{-6pt}
\begin{eqnarray*}
\hspace{18pt}\begin{array}{lcll}
E(A,B,\,c), & \mbox{with} & A\in\{a,\,c-a\}, & B\in\{b,\,c-b\},\\
E(A,B,2-c), & \mbox{with} & A\in\{1-a,1+a-c\}, & B\in\{1-b,1+b-c\},\\
E(A,B,1\!+\!a\!+\!b\!-\!c), & \mbox{with} & A\in\{a,1+b-c\}, & B\in\{b,1+a-c\},\\
E(A,B,1\!+\!c\!-\!a\!-\!b), & \mbox{with} & A\in\{1-a,c-b\}, & B\in\{1-b,c-a\},\\
E(A,B,1+a-b), & \mbox{with} & A\in\{a,\,1-b\}, & B\in\{1+a-c,c-b\},\\
E(A,B,1+b-a), & \mbox{with} & A\in\{1-a,\,b\}, &
B\in\{c-a,1+b-c\}.\hspace{20pt}\qed
\end{array}
\end{eqnarray*}

%This set contains bases of local solutions at all three singular points.
%\noindent %

The purpose of this paper is to present solutions (and relations between
them) of degenerate hypergeometric equations. These equations usually have
terminating or undefined hypergeometric solutions. In particular, if the
local exponent difference at a singular point is an integer, then that point
is either logarithmic or there is a terminating local solution at that
point. At a {\em logarithmic point} there is only one local solution of the
form $x^\lambda\!\left(1+\alpha_1x+\alpha_2x^2+\ldots\right)$, where $x$ is
a local parameter there. To get a basis of local solutions at a logarithmic
point, one has to use the function $\log x$.

Terminating solutions occur if $a$, $b$, $c-a$ or $c-b$ is an integer. Then
the hypergeometric equation has reducible monodromy group. %representation.
Conversely, if the monodromy group is reducible, the hypergeometric equation
has a solution which changes by a constant multiple under any monodromy
action. The logarithmic derivative of such a solution is a rational function
of $z$, which eventually means that the solution has a terminating series
expression. This is clear from the Kovacic algorithm in differential Galois
theory \cite{kovacic}, \cite[Section 4.3.4]{vdputsing}.

We recall briefly the role of the monodromy group. This group characterizes
analytic continuation of solutions of the hypergeometric equation along
paths in $\CC\setminus\{0,1\}$; see \cite[Section 3.9]{beukers}. Once a
basis of local solutions at a non-singular point is chosen, we get a
two-dimensional representation of the monodromy group. In general, the
monodromy group is generated by two elements, say, those corresponding to
paths circling $z=0$ or $z=1$ once. Fractional-linear transformations of
$E(a,b,c)$ do not change the monodromy group. We are especially interested
in the cases when a monodromy representation is a subgroup (up to
conjugation) of the following subgroups of $\mbox{GL}(2,\CC)$:
\[ %begin{eqnarray}
\GG_m = \left\{\left. \left(\begin{array}{cc} 1&0\\0&u \end{array}\right)
\right| u\in\CC^* \right\}, \qquad \GG_a = \left\{\left.
\left(\begin{array}{cc} 1&u\\0&1 \end{array}\right)\right| u\in\CC \right\}.
\] %end{eqnarray}
The groups $\GG_m$ and $\GG_a$ are isomorphic to the multiplicative group
$\CC^*$ and the additive group $\CC$, respectively. For Fuchsian equations,
the monodromy group is closely related to the differential Galois group
\cite[Theorem 5.3]{vdputsing}. The groups $\GG_m$, $\GG_a$ are examples of
possible differential Galois groups for hypergeometric equations.

Here is a convenient characterization of hypergeometric equations with
various types of monodromy groups.
\begin{theorem} \label{dendcases}
Consider a hypergeometric equation $E(a,b,c)$.
\begin{enumerate}
\item The monodromy group is irreducible if and only if none of the numbers
$a$, $b$, $c-a$, $c-b$ is an integer.%
\item If the monodromy group is abelian, then the sequence $a$, $b$, $c-a$,
$c-b$ contains at least two integers.%
\item Suppose that the sequence $a$, $1-b$, $c-a$, $1+b-c$ contains
precisely two integers. If those two integers are either both positive or
both non-positive, then a monodromy representation is a non-trivial subgroup
of $\GG_m$. Otherwise the monodromy group is not abelian.%
\item Suppose that the numbers $a$, $b$, $c$ are integers. If the sequence
$a$, $b$, $c-a$, \mbox{$c-b$} contains odd number of positive integers, then
the monodromy group is trivial. Otherwise a monodromy representation is a
non-trivial subgroup of $\GG_a$.%
\end{enumerate}
\end{theorem}
\proof The first statement is a direct consequence of \cite[Theorem
4.3.2]{gausspainleve}. Alternatively, this is \cite[Corollary
3.11]{beukers}. %A useful fact is \cite[Lemma 3.10]{beukers}

The second statement is \cite[Lemma 3.13]{beukers}. Alternatively, one can
go through cases (B.2), (B.2)$'$, (B.2)$''$, (C) of \cite[Theorem
4.3.3]{gausspainleve}. %In the setting of that theorem, we take
%$\rho_1=\sigma_1=0$, $\rho_2=1-c$, $\sigma_2=c-a-b$, $\tau_1=a$, $\tau_2=b$.

For the third statement, we may assume due to fractional-linear
transformations that $a,c-a$ are integers, and that $c\ge 1$. Then we are in
part (B) of \cite[Theorem 4.3.3]{gausspainleve}, with $\rho_1-\rho_2=c-1$
there. We have to check whether $0<\rho_1+\sigma_i+\tau_j\le \rho_1-\rho_2$
for those $i,j\in\{1,2\}$ with $\rho_1+\sigma_i+\tau_j\in\ZZ$. We have
$\rho_1+\sigma_i+\tau_j\in\{a,\,c-a,\,b,\,c-b\}$, therefore we are in case
(B.2) if the two integers are both positive, and we are in case (B.1) if
$a\le 0$ and $c-a>0$. The monodromy generators of case (B.1) do not commute,
and they are (up to conjugation) inside $\GG_m$ in case (B.2).
Fractional-linear transformations preserve the property that the two
integers in the sequence $a,1-b,c-a,1+b-c$ are either both positive, or both
non-positive.

For the last statement, we use part (C) of \cite[Theorem
4.3.3]{gausspainleve}. It is enough to show that if the sequence
$a,b,c-a,c-b$ contains odd number of positive integers then we are in case
(C.2). Up to fractional-linear transformations, we may assume that $a\le 0$
and that the other three integers are positive. Then $a\le 0<b<c$. We have
$\rho_1-\rho_2=c-1$, $\sigma_1-\sigma_2=c-a-b$, and for $i,j\in\{1,2\}$ we
have
\[
\rho_1+\sigma_i+\tau_j\in\{a,\,c-a,\,b,\,c-b\},\qquad
\rho_i+\sigma_1+\tau_j\in\{1-a,c-b,c-a,1-b\}.
\]
We check that $0<\rho_1+\sigma_i+\tau_j\le \rho_1-\rho_2$ for any
$i,j\in\{1,2\}$ in the context of part (B) in \cite{gausspainleve}, and
$0<\rho_i+\sigma_1+\tau_j\le \sigma_1-\sigma_2$ for any $i,j\in\{1,2\}$ in
the context of part (B)$'$. This means that we are indeed in case (C.2).
Fractional-linear transformations preserve the property that the sequence
$a,b,c-a,c-b$ contains odd number of positive integers. \qed

\begin{conclude} \label{abcases}
%Consider hypergeometric equation $E(a,b,c)$.
\begin{enumerate}
\item A monodromy representation of hypergeometric equation $E(a,b,c)$ is
(up to conjugation) a non-trivial subgroup of $\GG_m$ if and only if the
sequence $a,1-b,c-a,1+b-c$ contains precisely two integers, and those two
integers are either both positive or both non-positive.%
\item The monodromy group of $E(a,b,c)$ is trivial if and only if
$a,b,c\in\ZZ$ and the sequence $a$, $b$, $c-a$, $c-b$ contains odd number
of positive integers.%
\item The monodromy representation of $E(a,b,c)$ is (up to conjugation) a
non-trivial subgroup of $\GG_a$ if and only if $a,b,c\in\ZZ$ and the
sequence $a$, $b$, $c-a$, $c-b$ contains even number of positive integers.
\end{enumerate}
\end{conclude}
\proof Parts 2, 3, 4 of Theorem \ref{dendcases} describe mutually exclusive
cases for possible abelian monodromy groups of $E(a,b,c)$.\qed

In Sections \ref{terminating} through \ref{additive} below we study various
degenerate cases of Gauss hypergeometric functions. We concentrate on
relations between their hypergeometric and logarithmic solutions. Here we
present general forms of hypergeometric equations for each degeneration
type.
\begin{theorem} \label{dencases}
%Assume that $n,m,\ell\in\ZZnn$, %are non-negative integers,
%and assume that $a,b,c$ are complex numbers such that $a,b,c,c-a\not\in\ZZ$.
\begin{enumerate}
\item Suppose that a hypergeometric equation has terminating hypergeometric
solutions, non-abelian monodromy, and does not have logarithmic points. Up
to fractional-linear transformations, the hypergeometric equation has the
form $E(-n,a,c)$, where $n\in\ZZnn$ and $a,c,c-a\not\in\ZZ$.%
\item Suppose that a hypergeometric equation has logarithmic points, but
does not have terminating hypergeometric solutions. Up to fractional-linear
transformations, the hypergeometric equation has the form $E(a,b,m+1)$,
with $m\in\ZZnn$ and $a,b\not\in\ZZ$.%
\item Suppose that a hypergeometric equation has terminating solutions and
logarithmic points, and that the monodromy group is not abelian. Up to
fractional-linear transformations, the hypergeometric equation has the form
$E(a,-n,m+1)$, with $n,m\in\ZZnn$ and $a\not\in\ZZ$.%
\item Suppose that a monodromy representation of a hypergeometric equation
is a non-trivial subgroup of $\GG_m$. Up to fractional-linear
transformations, the hypergeometric equation has the form
$E(-n,a-m,-n-m)$, with $n,m\in\ZZnn$ and $a\not\in\ZZ$.%
\item Up to fractional-linear transformations, a hypergeometric equation
with the trivial monodromy group has the form $E(-n,\ell+1,-n-m)$, with
$n,m,\ell\in\ZZnn$.%
\item Suppose that a monodromy representation of a hypergeometric equation
is a non-trivial subgroup of $\GG_a$. Up to fractional-linear
transformations, the hypergeometric equation has the form
$E(-\ell,-n-\ell,-n-m-2\ell)$, with $n,m,\ell\in\ZZnn$.%
\end{enumerate}
\end{theorem}
\proof Consider a hypergeometric equation $E(a,b,c)$ with no restrictions on
the parameters $a,b,c$. In the first case, we may assume $b=-n\in\ZZnp$.
Since the monodromy group is not abelian and there are no logarithmic
points, the local exponent differences $1-c, a+n, c-a+n$ are not integers.

In the second case, we may choose $z=0$ as a logarithmic point. Then $c$ is
an integer, and we can choose it to be positive. Since there are no
terminating solutions, $a$ and $b$ are not integers by part 1 of Theorem
\ref{dendcases}.

In the third case, we choose $z=0$ as a logarithmic point as well, so $c$ is
a positive integer. Terminating solutions occur if $a$ or $b$ is an integer.
Assume that $b$ is an integer. Then $a\not\in\ZZ$ by part 4 of Theorem
\ref{dendcases}, and either $b\le 0$ or $c-b\le 0$ by part 3 of the same
theorem. While keeping $c$ positive, we can permute $b$ and $c-b$ by part 2
of Lemma \ref{frlints}, so we may assume that $b\le 0$.

In the fourth case, we may assume $a=-n$ and use part 1 of Corollary
\ref{abcases}, so one of the numbers $1-b,c-a,1+b-c$ is a non-positive
integer $-m$. Due to fractional-linear transformations we may assume that
$c-a=-m$. We are allowed to rename $b$ to $a-m$, for the purpose of
symmetric presentation in Section \ref{abelian}.

In the fifth case, the numbers $a,b,c$ are integers by part 2 of Corollary
\ref{abcases}. One of the numbers in the sequence $a,b,c-a,c-b$ has
different positivity than the others. Up to fractional-linear
transformations, we may assume that $b=\ell+1$ is a positive integer, and
$a=-n$, $c-a=-m$ are non-positive integers.

In the last case, the numbers $a,b,c$ are integers by part 3 of Corollary
\ref{abcases}. The four numbers in part 2 of Lemma \ref{frlints} cannot be
all negative, since their sum is zero. We may assume that $a,b,c-a,c-b$ are
all non-negative, and that $c-a\le b\le a\le 0$. Therefore we may set
$a=-\ell$, $b-a=-n$ and $c-a-b=-m$.\qed

\begin{table}
\begin{center} \begin{tabular}{|c|c|c|c|}
\hline & Kummer's & Terminating & Non-terminating \\
Case  & series & solutions & solutions \\
 \hline\hline
%0 & 24 & --- & 4, 4, 4, 4, 4, 4 \\
1 & 24 & 6+6 & 4, 4, 4 \\
2 & \ 12, 16 or 20\ & --- & 4, 4, 4 (and possibly 4, 4) \\
& 10, 13 or 16 & --- & 3, 3, 4 (and possibly 3, 3) \\
& 6, 8 or 10 & --- & 2, 2, 2; or 2, 3, 3; or 2, 2, 3, 3\\
3 & 16 or 20 & 6+4 or 8+4 & 3, 3 or 4, 4 \\
4 & 24 & 6+4,\ 6+4 & 4 \\
5 & 24 & \ 6+2,\ 6+2,\ 6+2\ & --- \\
6 & 10, 13 or 16 & 6+2; or 8+2; or 10+2 & 2; or 3; or 4 \\ \hline
\end{tabular} \vspace{-2pt}\end{center}
\caption{Kummer's solutions in degenerate cases} \label{figtab}
\end{table}
Table \ref{figtab} describes concisely the set of Kummer's solutions in each
degenerate case. The case numbers refer to Theorem \ref{dencases}. The
second column gives the total number of distinct well-defined Kummer's
series. In the third column, each terminating solution is represented by an
additive expression, where the first integer gives the number of terminating
hypergeometric expressions for the solution, and the second integer gives
the number of non-terminating expressions for the same solution. The last
column specifies solutions which have only non-terminating hypergeometric
expressions; usually these solutions have 4 hypergeometric expressions due
to Euler-Pfaff transformations. Multiple subcases are commented promptly
below in the following Section. %right after Lemma \ref{f21f23} here below.
We consider only relations between the 24 Kummer's solutions, so we do not
take into account quadratic or higher degree transformations, nor we
consider artificial identities with terminating series like
\[
\hpg21{a,-1}{b}{\,z}=\hpg21{c,-1}{d}{\frac{ad}{bc}z}.
\]
In particular, we consider two constant terminating $\hpgo21$ series (that
is, series with a zero upper parameter) as distinct if other parameters or
the argument are not equal. Correctness of Table \ref{figtab} is evident
from detailed considerations in Sections \ref{terminating} through
\ref{additive}.

\section{Some explicit facts}

Here we present some formulas which are useful in the following Sections. In
particular, we discuss %Lemmas \ref{f21f23} and \ref{f21ff} summarize
degenerations of Euler-Pfaff transformations, and introduce an alternative
normalization of Gauss hypergeometric function. %, see formula (\ref{althpgseries}).

Recall that Euler-Pfaff transformations \cite[Theorem 2.2.5]{specfaar}
identify the following hypergeometric series:
\begin{eqnarray} \label{flinear1}
\hpg{2}{1}{a,\,b\,}{c}{\,z} & = &
(1-z)^{c-a-b}\;\hpg{2}{1}{c-a,\,c-b\,}{c}{\,z\,}\\
\label{flinear2} & = & (1-z)^{-a}\;
\hpg{2}{1}{a,\,c-b\,}{c}{\frac{z}{z-1}} \\
\label{flinear3} &=& (1-z)^{-b}\; \hpg{2}{1}{c-a,\,b\,}{c}{\frac{z}{z-1}}.
\end{eqnarray}
These transformations hold when $c\not\in\ZZnp$. If the local exponent
differences $c-a-b$ and $a-b$ are non-zero, then these four series are
distinct. Note that the four series are local series at $z=0$. Recall that
Kummer's series for $E(a,b,c)$ with the argument $1-z$ or $1-1/z$ are local
series at $z=1$; Kummer's series with the argument $1/z$ or $1/(1-z)$ are
local series at $z=\infty$.

It turns out that if none of four series (\ref{flinear1})--(\ref{flinear3})
is terminating, which means
$a,b$, $c-a,c-b\not\in\ZZnp$, then no Kummer's series at % the other points
$z=1$, $z=\infty$ are equal to the left-hand of (\ref{flinear1}). To see
this, one may consider a solution basis and connection formulas in the
general case \cite[Section 2.9]{bateman} and in case 2 of Table \ref{figtab}
(see Section \ref{logarithmic} here).

We formulate a general conclusion as follows. Let $F$ denote a Gauss
hypergeometric function with the argument $\varphi(z)$ as in
(\ref{p1automor}), so that $\varphi(z)$ is a local parameter at a point
$P\in\{0,1,\infty\}$. Let $n(F)\in\{0,1,2\}$ denote the number of points in
the set $\{0,1,\infty\}\setminus\{P\}$ where the local exponent difference
for the corresponding hypergeometric equation is equal to zero. Then we have
exactly $4-n(F)$ distinct hypergeometric series expressions for $F$, unless
$F$ can be expressed as a terminating
series. %Terminating solutions have less hypergeometric expressions when
%there are points with the zero local exponent difference.

Now multiple cases of Table \ref{figtab} can be better clarified. Recall
that points with the zero local exponent difference are always logarithmic.
In case 2 we can have at most three logarithmic points, so we may have no
singular points with the zero local exponent difference (the first line
there), one such a point (the second line), or 2 or 3 such points (the third
line). In case 3 we have one logarithmic point where the local exponent
difference can be zero. In case 6 we have two logarithmic points.\\

On a few occasions we use the following normalization of Gauss
hypergeometric series:
\begin{equation} \label{althpgseries}
\hppg{a,\,b\,}{c}{\,z}:=
\sum_{k=0}^{\infty}\frac{\Gamma(a+k)\,\Gamma(b+k)}{\Gamma(c+k)\,\Gamma(1+k)}.
%=\frac{\Gamma(a)\,\Gamma(b)}{\Gamma(c)}\;\hpg{2}{1}{a,\,b\,}{c}{\,z}.
\end{equation}
This series is well defined when $a,b,c\not\in\ZZnp$; in this case
\[
\hppg{a,\,b\,}{c}{\,z}=\frac{\Gamma(a)\,\Gamma(b)}{\Gamma(c)}\;\hpg{2}{1}{a,\,b\,}{c}{\,z}.
\]
If $c=-N\in\ZZnp$ and $a,b\not\in\ZZnp$, then we should take limits of the
summands with singular gamma values:
\[
\hppg{a,\,b\,}{-N}{\,z}=\frac{\Gamma(a\!+\!N\!+\!1)\,\Gamma(b+\!N\!+\!1)}{(N+1)!}
\;\hpg{2}{1}{a\!+\!N\!+\!1,\,b+\!N\!+\!1\,}{N+2}{\,z}.
\]
If one of the parameters $a,b$ is a non-positive integer, but $c$ is a
smaller or equal integer, then the function in (\ref{althpgseries}) is also
well-defined if we agree to evaluate quotients of singular gamma values by
taking residues there. Then we have the following formulas.
\begin{lemma}\label{f21ff}
Suppose that $n,N\in\ZZnn$, and that $n\le N$. Then
\begin{eqnarray} \label{althpgrel}
\hppg{-n,a}{-N}{z}&=&(-1)^{N-n}\,\frac{\Gamma(a)\,N!}{n!}\,\hpg21{-n,a}{-N}{z}
\nonumber\\
&&+\frac{\Gamma(a\!+\!N\!+\!1)\,(N\!-\!n)!}{(N+1)!}\,z^{N+1}\,\hpg21{N-n+1,
a+N+1}{N+2}{z}.
\end{eqnarray}
Correct versions of Euler-Pfaff transformations are:
\begin{eqnarray} \label{eulera}
\hpg{2}{1}{-n,a}{-N}{\,z}&=&
(1-z)^{-a+n-N}\,\frac{(-1)^n\,(N-n)!}{N!\,\Gamma(-a-N)}\,\hppg{n-N,-a-N}{-N}{\,z}\\
&=&(1-z)^{-a}\,\frac{(-1)^n\,(N-n)!}{N!\,\Gamma(a)}\,\hppg{n-N,a}{-N}{\frac{z}{z-1}}\\
&=&(1-z)^{n}\,\hpg21{-n,-a-N}{-N}{\frac{z}{z-1}}.
\end{eqnarray}
\end{lemma}
\proof Formula (\ref{althpgrel}) is straightforward. For formula
(\ref{eulera}), we put $b=-n$, $c=-\nu-n$ in general Euler's formula
(\ref{flinear1}) and take the limit $\nu\to N-n$. For the other two
formulas, we take the same specialization and the same limit in Pfaff's
formulas (\ref{flinear2})-(\ref{flinear3}). \qed

For further convenience, we set forth the following two
functions:%\vspace{-3pt}
\begin{equation} \label{localat0}
H_1=\hppg{a,\,b\,}{c}{\,z},\qquad H_2=z^{1-c}\;\hppg{1+a-c,\,1+b-c}{2-c}{z}.
\end{equation}
In general, they form a basis of solutions for $E(a,b,c)$. Connection
formulas for other hypergeometric solutions of $E(a,b,c)$ can be
written as: %\vspace{-3pt}% Then the following formulas hold:
\begin{eqnarray} \label{altform1}
\hppg{a,\;b}{\!1+a+b-c}{1\!-\!z}&\equal&
\frac{H_1-H_2}{K},\\
\label{altform2}\hspace{-11pt}(1\!-\!z)^{c-a-b}\,
\hppg{c\!-\!a,\,c\!-\!b}{\!1\!+\!c\!-\!a\!-\!b}{1\!-\!z}&
\equal&\frac{1}{K}%{\Gamma(1+a-c)\,\Gamma(1+b-c)}
\left( H_1\,\frac{\sin\pi a\,
\sin\pi b}{\sin\pi(c-a)\sin\pi(c-b)}-H_2\right),\\
(-z)^{-a}\,\hppg{a,\,1+a-c}{1+a-b}{\,\frac1{z}\,}&\equal&
 \frac1{\sin\pi c}\left(H_1\,{\sin\pi b}+H_2\,{\bf e}^{i\pi c}\sin\pi(c-b)\right),\\
\label{altform4}(-z)^{-b}\,\hppg{1+b-c,\,b}{1+b-a}{\,\frac1{z}\,}&\equal&
 \frac1{\sin\pi c}\left(H_1\,{\sin\pi a}+H_2\,{\bf e}^{i\pi
 c}\sin\pi(c-a)\right),
\end{eqnarray}
where %\vspace{-4pt}
\[
K=\frac{\sin\pi c}{\pi}\,\Gamma(1+a-c)\,\Gamma(1+b-c).
\]
These formulas hold for analytic continuations of the $\bf F$-functions onto
the upper half-plane, like in \cite[Section 2.9]{bateman}.
%\cite[Theorem 2.3.2]{specfaar}, \cite[Section 15.3]{abrostegun} or \cite[Section 5.3]{nicoboek}.

We shall use the following consequences of the reflection formula
\cite[Theorem 1.2.1]{specfaar} for the gamma function:%\vspace{-4pt}
\begin{eqnarray} \label{deuler1}
\psi(x)-\psi(1-x)&\equal&-\frac{\pi}{\tan\pi x},\\ \label{deuler2}
\frac{\Gamma'(x)}{\Gamma(x)^2}&\equal&\frac{\Gamma'(1-x)}{\Gamma(x)\,\Gamma(1-x)}-\cos\pi
x\;\Gamma(1-x).
\end{eqnarray}
Recall that $\psi(x)=\Gamma'(x)/\Gamma(x)$.
\begin{lemma} \label{psitransf}
Suppose that $b\not\in\ZZ$ and $c\not\in\ZZnp$. Then for $|z|<1$ (and
eventually, after analytic continuation) we have
\[
\sum_{k=0}^{\infty} \frac{(a)_k(b)_k}{(c)_k\,k!}\, \psi(b+k)\,z^k=
\sum_{k=0}^{\infty} \frac{(a)_k(b)_k}{(c)_k\,k!}\,\psi(1\!-\!b\!-\!k)\,z^k-
\frac{\pi}{\tan\pi b}\,\hpg21{a,b}{c}{z}.\hspace{-4pt}
\]
\end{lemma}
\proof The formula follows by applying (\ref{deuler1}) termwise. \qed

\section{Terminating hypergeometric series}
\label{terminating}

Here we consider hypergeometric equations which have terminating
hypergeometric solutions and non-abelian monodromy group, but do not have
logarithmic points. We are in case 1 of Theorem \ref{dencases}. A general
equation is $E(-n,a,c)$, where $n$ is a non-negative integer, and $a,c,c-a$
are not integers. All 24 Kummer's solutions are well defined. The monodromy
group of this equation is reducible, because a terminating hypergeometric
solution spans an invariant subspace.

It turns out that there are terminating hypergeometric solutions at each
singular point of $E(-n,a,c)$. %the hypergeometric equation.
All terminating solutions lie in the one-dimensional invariant subspace.
Here is their identification: %\vspace{-1pt} %they are proportional:
\begin{eqnarray} \label{termseries1}
\hpg{2}{1}{-n,\,a}{c}{\,z\,} &=&(1-z)^n\;\hpg{2}{1}{-n,\,c-a}{c}{\frac{z}{z-1}}\\
&=&\frac{(a)_n}{(c)_n}\,(-z)^n\,\hpg{2}{1}{-n,\;1-n-c}{1-n-a}{\,\frac1{z}\,}\\
&=&\frac{(a)_n}{(c)_n}\,(1-z)^n\,\hpg{2}{1}{-n,\;c-a}{1-n-a}{\frac1{1-z}}\\
&=&\frac{(c-a)_n}{(c)_n}\,z^n\,\hpg{2}{1}{-n,\;1-n-c}{1-n+a-c}{\,1-\frac1{z}\,}\\
\label{termseries6}&=&\frac{(c-a)_n}{(c)_n}\;\hpg{2}{1}{-n,\;a}{1-n+a-c}{\,1-z\,}.
\end{eqnarray}
These formulas can be proved using the following two transformations a few
times: rewriting a terminating hypergeometric sum in the opposite direction,
and Pfaff's formula (\ref{flinear2}). Application of Euler's formula
(\ref{flinear1}) to the above series gives non-terminating hypergeometric
expressions for the same function. We have 6 terminating and 6
non-terminating hypergeometric expressions for this solution.

The six expressions (\ref{termseries1})--(\ref{termseries6}) are valid for
any terminating Gauss hypergeometric series, if only the numbers $a$, $c$,
$c-a$ are not integers in the interval $[1-n,0]$. Any terminating series can
be interpreted as an isolated Jacobi, Meixner or Meixner-Pollaczek
polynomial \cite[Sections 1.7,1.8,1.9]{koekswart}:
\begin{eqnarray}
\hpg{2}{1}{-n,\,a}{c}{\,z\,}&=&\frac{n!}{(c)_n}\;P^{(c-1,a-c-n)}_n(1-2z)\\
&=&M_n\!\left(-a;\, c, \frac1{1-z}\right)\\
&=&\frac{n!}{(c)_n(1\!-\!z)^{n/2}}\;P^{(c/2)}_{n}\!\left({\frac{ic}{2}-ia;\,
\frac{i}{2}\log(1-z)}\right)
\end{eqnarray}
We get back to the specified equation $E(-n,a,c)$. Its
non-terminating solutions are%\vspace{-2pt}
\begin{eqnarray*}
z^{1-c}\,(1\!-\!z)^{c-a+n}\,\hpg{2}{1}{1+n,1-a}{2-c}{z},\quad
z^{1-c}\,(1\!-\!z)^{c-a+n}\,\hpg{2}{1}{1+n,\,1-a}{1+n+c-a}{1\!-\!z},\\
(-z)^{-c-n}\,(1-z)^{c-a+n}\,\hpg{2}{1}{1+n,\,c+n\,}{1+n+a}{\,\frac1z\,}.
\end{eqnarray*}
For each of these functions, there are 4 non-terminating hypergeometric
expressions by Euler-Pfaff transformations. This exhausts the 24 Kummer's
solutions. Connection relations are evident from the following
formulas:%\vspace{-1pt}
\begin{eqnarray*}
\hpg{2}{1}{1+n,\,1-a}{1+n+c-a}{1-z}&\equal&\frac{(c-a)_{n+1}}{(c-1)_{n+1}}\,
\hpg{2}{1}{1+n,1-a}{2-c}{z}\\ &&\hspace{-48pt}+\frac{\Gamma(1+n+c-a)\,
\Gamma(1-c)}{\Gamma(1-a)\;n!}\,z^{c-1}\,(1-z)^{a-c-n}\hpg{2}{1}{-n,\,a}{c}{z},\\
\hpg{2}{1}{1+n,\,c+n\,}{1+n+a}{\,\frac1z\,}&\equal&
\frac{(a)_{n+1}}{(c-1)_{n+1}}\,(-z)^{n+1}\,\hpg{2}{1}{1+n,1-a}{2-c}{z}\\
&&\hspace{-48pt}+\frac{\Gamma(1+n+a)\,\Gamma(1-c)}{\Gamma(1+a-c)\;n!}
\,(-z)^{c+n}\,(1-z)^{a-c+n}\,\hpg{2}{1}{-n,\,a}{c}{z}.
\end{eqnarray*}

\section{General logarithmic solutions}
\label{logarithmic}

Here we consider hypergeometric equations which have logarithmic points, but
do not have terminating hypergeometric solutions. We are in case 2 of
Theorem \ref{dencases}. A general hypergeometric equation of this kind is
$E(a,b,m+1)$, where $m$ is a non-negative integer, and $a,b$ are not
integers.

The functions $H_1$, $H_2$ in (\ref{localat0}) coincide in this case. The
corresponding $\hpgo21$ series either coincide (if $m=0$), or only one of
them is well-defined (if $m\ge 1$).
%There is only one Gauss hypergeometric series
Formulas for a second independent local solution at $z=0$ are not pretty,
but it must be important to have them. We choose to identify logarithmic
solutions with
\begin{equation} \label{logsol2}
U_1=(-1)^{m+1}\,m!\,\frac{\Gamma(a-m)\,\Gamma(b-m)}{\Gamma(a+b-m)}\,
\hpg{2}{1}{a,\;b}{\!a+b-m}{1\!-\!z}.
\end{equation}
%Recall that $\psi(x)=\Gamma'(x)/\Gamma(x)$.
\begin{theorem} \label{loggent}
The function $U_1$ has the following expressions:
\begin{eqnarray} \label{logsol3}
U_1&\equal&(-1)^{m+1}m!\,\frac{\Gamma(1-a)\,\Gamma(1-b)}{\Gamma(m+2-a-b)}
(1\!-\!z)^{m+1-a-b}\,\hpg{2}{1}{m\!+\!1\!-\!a,\,m\!+\!1\!-\!b}{\!m+2-a-b}{1\!-\!z}\nonumber\\
&&-\frac{\pi\,\sin\pi(a+b)}{\sin\pi a\,\sin\pi b}\,\hpg{2}{1}{a,\,b}{m+1}{\,z},\\
\label{u1ab}&\equal&(-1)^{m+1}m!\,\frac{\Gamma(a-m)\,\Gamma(1-b)}{\Gamma(1+a-b)}
(-z)^{-a}\hpg{2}{1}{a,\,a-m}{\!1+a-b}{\frac1z}\nonumber\\
&&-\frac{\pi\,{\bf e}^{-i\pi b}}{\sin\pi b}\,\hpg{2}{1}{a,\,b}{m+1}{\,z}\\
\label{logsol1} &\equal& \hpg{2}{1}{a,b}{m+1}{\,z}\log z
+\frac{(-1)^{m+1}m!\,(m\!-\!1)!}{(1-a)_m\,(1-b)_m}\,z^{-m}
\,\sum_{k=0}^{m-1}\frac{(a\!-\!m)_k(b\!-\!m)_k}{(1-m)_k\;k!}{z^k}\nonumber\\
&&+\sum_{k=0}^{\infty}\frac{(a)_k(b)_k}{(m+1)_k\,k!}\big(
\psi(a\!+\!k)+\psi(b\!+\!k)-\psi(m\!+\!k\!+\!1)-\psi(k\!+\!1)\big)\,{z^k}\\
\label{logsol3a} &\equal&\hpg{2}{1}{a,\,b}{m+1}{\,z}\log z
-\frac{\pi\sin\pi(a+b)}{\sin\pi a\,\sin\pi b}\,\hpg{2}{1}{a,\,b}{m+1}{\,z}\nonumber\\
&&+\frac{(-1)^{m+1}m!\,(m-1)!}{(1-a)_m\,(1-b)_m}\,z^{-m}\,(1\!-\!z)^{m+1-a-b}\,
\,\sum_{k=0}^{m-1}\frac{(1-a)_k(1-b)_k}{(1-m)_k\;k!}{z^k}\nonumber\\
\nonumber &&+(1\!-\!z)^{m+1-a-b}\sum_{k=0}^{\infty}\!
\frac{(m\!+\!1\!-\!a)_k(m\!+\!1\!-\!b)_k}{(m+1)_k\,k!}\times\\
&&\hspace{24pt}\big(\psi(m\!+\!k\!+\!1\!-\!a)+\psi(m\!+\!k\!+\!1\!-\!b)
-\psi(m\!+\!k\!+\!1)-\psi(k\!+\!1)\big)\,{z^k}\\
&\equal&\hpg{2}{1}{a,\,b}{m+1}{\,z}\log\frac{z}{1-z}
-\frac{\pi}{\tan\pi b}\,\hpg{2}{1}{a,\,b}{m+1}{\,z}\nonumber\\
\label{u1ab2}&&+\frac{(-1)^{m+1}m!(m-1)!}{(1-a)_m(1-b)_m}\,z^{-m}(1\!-\!z)^{m-a}
\sum_{k=0}^{m-1}\frac{(a\!-\!m)_k(1\!-\!b)_k}{(1-m)_k\;k!}{\frac{z^k}{(z\!-\!1)^k}}\nonumber\\
&&+(1\!-\!z)^{-a}\sum_{k=0}^{\infty}\frac{(a)_k(m\!+\!1\!-\!b)_k}{(m+1)_k\,k!}\times\nonumber\\
&&\hspace{24pt}\big(\psi(a\!+\!k)+\psi(m\!+\!k\!+\!1\!-\!b)-\psi(m\!+\!k\!+\!1)-\psi(k\!+\!1)\big)
{\frac{z^k}{(z\!-\!1)^k}}.
\end{eqnarray}
\end{theorem}
\proof Formulas (\ref{logsol3}) and (\ref{u1ab}) are special cases of
connection formulas 2.9.(33) and 2.9.(25) in \cite{bateman}, respectively.

To prove formula (\ref{logsol1}), we apply formula (\ref{altform1}) to the
equation $E(a,b,m+1)$ and use expression (\ref{logsol2}). We take the limit
$c\to m+1$ on the right-hand of (\ref{altform1}) by l'Hospital's rule and
arrive at
\begin{eqnarray} \label{askroy}
U_1&\equal&-\frac{m!}{\Gamma(a)\,\Gamma(b)}\,\sum_{k=0}^{\infty}
\frac{\Gamma(a+k)\,\Gamma(b+k)\,\Gamma'(m+k+1)}{k!\;\Gamma(m+k+1)^2}\,z^k\nonumber \\
&&+\frac{m!\;z^{-m}\,\log z}{\Gamma(a)\,\Gamma(b)}\,\sum_{k=0}^{\infty}
\frac{\Gamma(a-m+k)\,\Gamma(b-m+k)}{k!\;\Gamma(1-m+k)}\,z^k\nonumber \\
&&+\frac{m!\;z^{-m}}{\Gamma(a)\Gamma(b)}\sum_{k=0}^{\infty}
\frac{\Gamma'(a\!-\!m\!+\!k)\,\Gamma\!(b\!-\!m\!+\!k)+
\Gamma\!(a\!-\!m\!+\!k)\,\Gamma'(b\!-\!m\!+\!k)}{k!\;\Gamma(1-m+k)}\,z^k\nonumber \\
&&-\frac{m!\;z^{-m}}{\Gamma(a)\,\Gamma(b)}\lim_{c\to
m+1}\,\sum_{k=0}^{\infty}
\frac{\Gamma(1\!+\!a\!-\!c\!+\!k)\,\Gamma(1\!+\!b\!-\!c\!+\!k)\,\Gamma'(2\!-\!c\!+\!k)}
{k!\;\Gamma(2-c+k)^2}\,z^k.
\end{eqnarray}
The first $m$ terms of the second and of the third series are zero. The
second series becomes the first term of (\ref{logsol1}). To compute the
first $m$ terms of the fourth series, we set $x=2-c+k$ in formula
(\ref{deuler2}) and get
\[
\lim_{c\to m+1}
\frac{\Gamma'(2-c+k)}{\Gamma(2-c+k)^2}=(-1)^{m-k}\,\Gamma(m-k)
=(-1)^m\frac{(m-1)!}{(1-m)_k}.
\]
The first $m$ terms of the fourth series in (\ref{askroy}) form the second
term of formula (\ref{logsol1}). The last term in (\ref{logsol1}) is
obtained by combining the remaining terms of the third and fourth series in
(\ref{askroy}), and all terms of the first series in (\ref{askroy}). This
proof of (\ref{logsol1}) follows closely the derivation in
\cite[pg.~82-84]{specfaar}; there are a few misprints in formulas there.

To prove (\ref{logsol3a}), we express the first term in %on the right-hand side of
(\ref{logsol3}) as a logarithmic function using (\ref{logsol1}). Then we
simplify the term with $\log z$ by Euler's formula.

To prove (\ref{u1ab2}), we apply Pfaff's transformation to the first term of
(\ref{u1ab}) and substitution $z\mapsto z/(z-1)$ in (\ref{logsol2}) and
(\ref{logsol1}). Then we simplify the logarithmic term by Pfaff's formula,
observing that (for $z\in\CC$ in the upper half-plane)%\vspace{-2pt}
\[
\hspace{68pt} \log\frac{z}{z-1}=\log\frac{z}{1-z}-i\pi,\qquad
\frac{\pi\,{\bf e}^{-i\pi b}}{\sin\pi b} =\frac{\pi}{\tan\pi b}-i\pi.%
 \hspace{68pt}\qed
\]
\\

More expressions for $U_1$ can be obtained by interchanging $a$, $b$ in
(\ref{u1ab}) and (\ref{u1ab2}), by applying Euler-Pfaff transformations to
individual hypergeometric functions in
\mbox{(\ref{logsol2})--(\ref{u1ab2})}, and by applying Lemma \ref{psitransf}
to sums with the $\psi$-function. In this way logarithmic solutions at $z=0$
can be related to any other well-defined hypergeometric series at $z=1$ or
$z=\infty$. One may check that well-defined hypergeometric series at
different points are independent.

The points $z=1$ and $z=\infty$ are not logarithmic if and only if
$a+b,\,a-b\not\in\ZZ$. We have either 4 undefined Kummer's series at $z=0$
(if \mbox{$m>0$}) or 4 pairs of coinciding hypergeometric series there (if
$m=0$). In the latter case, Euler's transformation (\ref{flinear1}) acts
trivially on some hypergeometric series at $z=1$ and $z=\infty$, for example
on (\ref{logsol2}). % for $U_1$.
Then we have 3 (rather than usual 4) distinct hypergeometric expressions for
each Gauss hypergeometric function representable by well-defined series at
$z=1$ and $z=\infty$. In any case (when $z=1$, $z=\infty$ are not
logarithmic), we have 5 different Gauss hypergeometric solutions represented
by 16 or 20 distinct Kummer's series.

If either $a+b$ or $a-b$ is an integer, then the equation $E(a,b,m+1)$ has
one other logarithmic point. For example, if $\ell\in\ZZnn$ and $a=b+\ell$
then the point $z=\infty$ is logarithmic. A power series solution there is
\[
z^{-a}\,\hppg{a,\,a-m}{\ell+1}{\,\frac1z}=z^{-b}\,\hppg{b,\,b-m}{1-\ell}{\,\frac1z}
\]
Just as we obtained solution (\ref{logsol1}) from coinciding solutions
$H_1=H_2$ in (\ref{althpgseries}) with $c=m+1$, we have the following local
solution at $z=\infty$:
\begin{eqnarray}
z^{-a}\,\hpg{2}{1}{a,a-m}{\ell+1}{\frac1z}\log\frac1z
+\frac{(-1)^{\ell+1}\ell!\,(\ell\!-\!1)!}{(1-a)_\ell\,(m+\!1-\!a)_\ell}\,z^{-b}
\,\sum_{k=0}^{\ell-1}\frac{(b)_k(b\!-\!m)_k}{(1-\ell)_k\;k!}z^{-k}\nonumber\\
+z^{-a}\sum_{k=0}^{\infty}\frac{(a)_k(a\!-\!m)_k}{(\ell+1)_k\,k!}\big(
\psi(a\!+\!k)+\psi(a\!-\!m\!+\!k)-\psi(\ell\!+\!k\!+\!1)-\psi(k\!+\!1)\big)z^{-k}.
\hspace{-4pt}
\end{eqnarray}
This function has to be identified with the following constant multiple of
$U_1$:
\[
(-1)^{\ell+1}\,\ell!\,\frac{\Gamma(b)\,\Gamma(b-m)}{\Gamma(a+b-m)}\,
\hpg{2}{1}{a,\;b}{\!a+b-m}{1\!-\!z}.
\]
To see this, use connection formula \cite[2.9.(36) with a misprint:
$\Gamma(c)$ must be replaced by $\Gamma(b)$]{bateman} before taking the
limit $b\to a-\ell$ in the analogue of (\ref{altform1}). If $m\neq 0$ and
$\ell\neq 0$, then 8 Kummer's series are undefined; other 16 Kummer's series
are distinct and represent 4 different Gauss hypergeometric functions. If
$m\neq 0$, $\ell=0$, then 4 Kummer's series at \mbox{$z=0$} are undefined,
and there are 4 coinciding pairs of them at \mbox{$z=\infty$}. Other 4
Kummer's series at $z=\infty$ are distinct expressions of one Gauss
hypergeometric function. There are 3 more Gauss hypergeometric solutions,
each represented by 3 distinct Kummer's series at $z=0$ or $z=1$. If $m=0$,
$\ell\neq 0$, we have a similar situation. If $m=\ell=0$, then there are 3
distinct Kummer's series at $z=0$ and at $z=\infty$, and 4 distinct series
at $z=1$. The series at $z=0$ or at $z=\infty$ represent a single Gauss
hypergeometric function; the series at $z=1$ represent two different
functions.

All three singular points of $E(a,b,m+1)$ are logarithmic if and only if $a$
and $b$ are rational numbers with the denominator 2. Then we have 3 distinct
Gauss hypergeometric functions, one for each singular point. The number of
their hypergeometric expressions depends on the presence of singular points
with the zero local exponent difference; recall discussion of Table
\ref{figtab}. If all three local exponent differences are 0, then each Gauss
hypergeometric solution is represented by just two distinct Kummer's series.
In this case the equation is $E(1/2,1/2,1)$; its hypergeometric solutions
are related to the well-known complete elliptic integral $K(k)$; see
\cite[(3.2.3)]{specfaar}. Other renowned complete elliptic integral is
$E(k)$ \cite[(3.2.14)]{specfaar}, which is expressible via solutions of
$E(-1/2,1/2,1)$. This equation has one local exponent difference equal to 0,
so the number of distinct Kummer's series is 10. An example of
hypergeometric equation with 13 distinct Kummer's series is
$E(-1/2,-1/2,1)$; it has two local exponent differences equal to 0.

\section{Logarithmic and terminating solutions}

Here we consider hypergeometric equations which have logarithmic points and
terminating hypergeometric solutions, but the monodromy group is
non-abelian. By part 3 of Theorem \ref{dencases},
%Hypergeometric equations with logarithmic and terminating solutions, and with
%abelian monodromy group are considered in Section \ref{additive}.
a general hypergeometric equation of this kind is $E(a,-n,m+1)$, where
$n,m\in\ZZnn$ and $a\not\in\ZZ$. We are in the non-abelian case of part 3 of
Theorem \ref{dendcases}. The point $z=0$ is logarithmic; the points $z=1$,
$z=\infty$ are not logarithmic.

The terminating solution $\hpg21{-n,\,a}{m+1}{z}$ has 6 terminating
expressions (\ref{termseries1})--(\ref{termseries6}) with $n+1$ terms, but
there are 2 extra terminating expressions with $m+n+1$ terms if $m\neq 0$:
\begin{eqnarray*} \label{ttermseries1}
\hpg{2}{1}{-n,\,a}{m+1}{\,z\,} &=&(1-z)^n\;\hpg{2}{1}{-n,\,m+1-a}{m+1}{\frac{z}{z-1}}\\
&=&\frac{m!\,(a)_n}{(m+n)!}\,(-z)^n\,\hpg{2}{1}{-n,\;-m-n}{1-n-a}{\,\frac1{z}\,}\\
&=&\frac{m!\,(a)_n}{(m+n)!}\,(1-z)^n\,\hpg{2}{1}{-n,\;m+1-a}{1-n-a}{\frac1{1-z}}\\
&=&\frac{m!\,(m\!+\!1\!-\!a)_n}{(m+n)!}\,z^n\,\hpg{2}{1}{-n,\;-m-n}{a-m-n}{\,1-\frac1{z}\,}\\
\label{ttermseries6}&=&\frac{m!\,(m\!+\!1\!-\!a)_n}{(m+n)!}\;\hpg{2}{1}{-n,\;a}{a-m-n}{\,1-z\,}\\
&=&\frac{m!\,(a)_n}{(m+n)!}\,(-z)^{-m}(1-z)^{m+n}\,\hpg{2}{1}{-m-n,\;1-a}{1-n-a}{\frac1{1-z}}\\
&=&\frac{m!\,(m\!+\!1\!-\!a)_n}{(m+n)!}\,z^{-m}\,\hpg{2}{1}{-m-n,\;a-m}{a-m-n}{\,1-z\,}.
\end{eqnarray*}
This solution has 4 non-terminating hypergeometric expressions with the
argument $z$, $z/(z-1)$, $1/z$ or $1-1/z$, due to Euler's formula
(\ref{flinear1}).  There are 2 other Gauss hypergeometric solutions,
represented by 6 (if $m>0$) or 8 (if $m=0$) distinct non-terminating
Kummer's series at $z=1$ and $z=\infty$. Like in Section \ref{logarithmic},
we miss 4 Kummer's series at $z=0$, which are either undefined  or coincide
with listed terminating series.

The logarithmic solution $U_1$ of Section \ref{logarithmic} is not defined
in this case, since some values of the $\psi$-function become infinite in
formulas (\ref{logsol2})--(\ref{u1ab2}).  We should either apply those
formulas to the equation $E(m+1-a,m+n+1,m+1)$, or consider the following
solution of $E(a,-n,m+1)$, well defined for $b=-n$ by expression
(\ref{logsol3a}):
%There appear to be no relations between terminating and logarithmic solutions.
\begin{equation} \label{u2exp0}
U_2:=U_1+\frac{\pi\sin\pi(a+b)}{\sin\pi a\,\sin\pi  b}
\,\hpg{2}{1}{a,\,b}{m+1}{\,z}
\end{equation}
The following Theorem presents various expressions for this function.
\begin{theorem} \label{logtermt}
The function $U_2$ with $b=-n$ has the following expressions:
\begin{eqnarray} \label{u2exp1}
U_2&\equal&\frac{(-1)^{m+1}\,m!\,n!}{(1-a)_{m+n+1}}
(1\!-\!z)^{m+n+1-a}\,\hpg{2}{1}{m\!+\!1\!-\!a,\,m\!+\!n\!+\!1}{\!m+n+2-a}{1\!-\!z}\\
\label{u2exp2}&\equal&\frac{(-1)^{m+1}\,m!\,n!}{(a-m)_{m+n+1}}
(-z)^{-a}\hpg{2}{1}{a,\,a-m}{a+n+1}{\,\frac1z}+\frac{\pi\,{\bf e}^{i\pi
a}}{\sin\pi a}\,\hpg{2}{1}{-n,\,a}{m+1}{\,z}\\
\label{u2exp3}&\equal&\hpg21{-n,\,a}{m+1}{z}\,\log z
+\frac{\pi}{\tan\pi a}\,\hpg21{-n,\,a}{m+1}{z}\nonumber\\
&&+\frac{(-1)^{m+1}\,m!\,n!}{(1-a)_m}\,z^{-m}\,\sum_{k=0}^{m-1}
\frac{(a-m)_k\,(m-k-1)!}{(m+n-k)!\,k!}\,z^k\nonumber\\ \mbox{\vspace{20pt}}
&&\hspace{-1pt}+\sum_{k=0}^n\!\frac{m!\,n!\,(a)_k}{(m\!+\!k)!(n\!-\!k)!\,k!}\!
\left(\psi(a\!+\!k)+\psi(n\!-\!k\!+\!1)-\psi(m\!+\!k\!+\!1)-\psi(k\!+\!1)\right)(-z)^k\nonumber\\
&&+(-1)^n\,n!\,m!\sum_{k=n+1}^{\infty}\!\frac{(a)_k(k-n-1)!}{(m+k)!\,k!}\,z^k\\
\label{u2exp4}&\equal&\hpg{2}{1}{-n,\,a}{m+1}{\,z}\log z\nonumber\\
&&+\frac{(-1)^{m+1}\,m!}{(1-a)_m(m\!+\!n)!}\,z^{-m}(1\!-\!z)^{m+n+1-a}\,
\sum_{k=0}^{m-1}\frac{(1-a)_k(m\!-\!k\!-\!1)!(n\!+\!k)!}{k!}(-z)^{k}\nonumber\\
&&+\frac{m!}{(m+n)!}\,(1-z)^{m+n+1-a}\,
\sum_{k=0}^{\infty}\frac{(m\!+\!1\!-\!a)_k(m\!+n\!+\!k)!}{(m+k)!\,k!}\times\nonumber\\
&&\hspace{24pt}\big(\psi(m\!+\!k\!+\!1\!-\!a)+\psi(m\!+\!n\!+\!k\!+\!1)
-\psi(m\!+\!k\!+\!1)-\psi(k\!+\!1)\big)\,z^k\\
\label{u2exp5}&\equal&\hpg{2}{1}{-n,\,a}{m+1}{\,z}\log
\frac{z}{1-z}+\frac{\pi}{\tan\pi a}
\,\hpg{2}{1}{-n,\,a}{m+1}{\,z}\nonumber\\
&&+\frac{(-1)^{m+1}\,m!}{(1-a)_m(m\!+\!n)!}\,z^{-m}(1\!-\!z)^{m-a}\,
\sum_{k=0}^{m-1}\frac{(a-m)_k(m\!-\!k\!-\!1)!(n\!+\!k)!}{k!}\frac{z^k}{(1-z)^k}\nonumber\\
&&+\frac{m!}{(m+n)!}\,(1-z)^{-a}\,
\sum_{k=0}^{\infty}\frac{(a)_k(m\!+n\!+\!k)!}{(m+k)!\,k!}\times\nonumber\\
&&\hspace{24pt}\big(\psi(a\!+\!k)+\psi(m\!+\!n\!+\!k\!+\!1)
-\psi(m\!+\!k\!+\!1)-\psi(k\!+\!1)\big)\,\frac{z^k}{(z-1)^k}\\
\label{u2exp6}&\equal&\hpg21{-n,\,a}{m+1}{z}\,\log \frac{z}{1-z}\nonumber\\
&&+\frac{(-1)^{m+1}\,m!\,n!}{(1-a)_m}\,z^{-m}\,(1-z)^{n+m}\sum_{k=0}^{m-1}
\frac{(1-a)_k\,(m-k-1)!}{(m+n-k)!\,k!}\,\frac{z^k}{(z-1)^k}\nonumber\\
&&\hspace{-1pt}+m!\,n!\,(1-z)^n\,\sum_{k=0}^n\!\frac{(m\!+\!1\!-\!a)_k}{(m\!+\!k)!(n\!-\!k)!\,k!}
\times\nonumber\vspace{-1pt}\\ &&\hspace{24pt}
\left(\psi(m\!+\!k\!+\!1\!-\!a)+\psi(n\!-\!k\!+\!1)-\psi(m\!+\!k\!+\!1)-\psi(k\!+\!1)\right)
\frac{z^k}{(1-z)^k}\nonumber\\
&&+m!\,n!\,(z-1)^n\,\sum_{k=n+1}^{\infty}\!\frac{(m\!+\!1-\!a)_k(k-n-1)!}{(m+k)!\,k!}
\,\frac{z^k}{(z-1)^k}.
\end{eqnarray} %\vspace{-8pt}
\end{theorem}
\proof To show the first two formulas we apply, respectively,
(\ref{logsol3}) or (\ref{u1ab}) to expression (\ref{u2exp0}). Then we
collect the two terms with $\hpg21{a,b}{m+1}z$ and take the limit $b\to -n$.

To prove (\ref{u2exp3}), we apply Lemma \ref{psitransf} to expression
(\ref{logsol1}) and arrive at%\vspace{-1pt}
 \begin{eqnarray*}
U_2&\equal& \hpg{2}{1}{a,\,b}{m+1}{\,z}\log z+\left(\frac{\pi\sin\pi(a+b)}
{\sin\pi a\,\sin\pi b}-\frac{\pi}{\tan\pi b}\right)\hpg21{a,\,b}{m+1}{z}\\
&&+\frac{(-1)^{m+1}m!\,(m\!-\!1)!}{(1-a)_m\,(1-b)_m}\,z^{-m}
\,\sum_{k=0}^{m-1}\frac{(a\!-\!m)_k(b\!-\!m)_k}{(1-m)_k\;k!}{z^k}\nonumber\\
&&+\sum_{k=0}^{\infty}\frac{(a)_k(b)_k}{(m+1)_k\,k!}\big(
\psi(a\!+\!k)+\psi(1\!-\!b\!-\!k)-\psi(m\!+\!k\!+\!1)-\psi(k\!+\!1)\big)\,{z^k}.
\end{eqnarray*}
Now we take the limit $b\to-n$ of each term. In particular, for $k\ge n+1$
we apply formula (\ref{deuler1}) with $x=b+k$, and then formula (\ref{deuler2})
to get %\vspace{-6pt}
\begin{eqnarray*}
\lim_{b\to-n} (b)_k\,\psi(1\!-\!b\!-\!k)&\equal& \lim_{b\to-n}
\frac{\Gamma(b+k)}{\Gamma(b)}\,\left(\psi'(b+k) %\frac{\Gamma'(b+k)}{\Gamma(b+k)}
+\frac{\pi\cos\pi b}{\sin\pi b}\right)\\ &\equal& \lim_{b\to-n}\left(
\frac{\Gamma'(b+k)}{\Gamma(b)}+\cos\pi b\;\Gamma(1-b)\,\Gamma(b+k)\right)\\
&\equal&(-1)^n\,n!\,(k-n-1)!.
\end{eqnarray*}
Note that the last sum in (\ref{u2exp3}) can be easily missed out; see
\cite[pg. 84]{specfaar}.

Formulas (\ref{u2exp4}), (\ref{u2exp5}) are direct consequences of and
(\ref{logsol3a}), (\ref{u1ab2}) applied to (\ref{u2exp0}). To get formula
(\ref{u2exp6}), consider (\ref{u1ab2}) with interchanged roles of the
parameters $a$ and $b$, so that $a$ approaches the integer $-n$ (and $b$ can
be renamed to $a$ at some point). Then we apply Lemma \ref{psitransf}
similarly as in the proof of (\ref{u2exp3}).\qed \vspace{-4pt}

\section{Completely reducible monodromy group}
\label{abelian}

Here we consider hypergeometric equations with completely reducible but
non-trivial monodromy group. Up to conjugation, the monodromy representation
is a subgroup of $\GG_m$. By part 4 of Theorem \ref{dencases}, a general
hypergeometric equation of this type is $E(-n,a-m,-n-m)$, where $n,m$ are
non-negative integers, and $a\not\in\ZZ$.

Since the monodromy group is completely reducible, there is a basis of
terminating solutions of $E(-n,a-m,-n-m)$. Such a basis is
\begin{equation} \label{abbasis}
\hpg{2}{1}{-n,\,a-m}{-n-m}{\,z\,}, \qquad
(1-z)^{-a}\;\hpg21{-m,-a-n}{-n-m}{\,z\,}.
\end{equation}
Although these two hypergeometric series seem to be equal by standard
Euler's formula (\ref{flinear1}), the correct Euler's transformation in this
situation is formula (\ref{eulera}) of Lemma \ref{f21ff}. Alternative
terminating expressions of the basis solutions (\ref{abbasis}) are obtained
by using formulas (\ref{termseries1})--(\ref{termseries6}). For example,
\begin{eqnarray*}
\hpg{2}{1}{-n,a-m}{-n-m}{\,z\,} &=& (1-z)^n\;\hpg{2}{1}{-n,-a-n}{-n-m}{\frac{z}{z-1}}\\
&=&\frac{m!\,(a-m)_n}{(n+m)!}\,z^n\,\hpg{2}{1}{-n,\;m+1}{\!1-a+m-n}{\,\frac1{z}\,}\\
&=&\frac{m!\,(a-m)_n}{(n+m)!}\,(z-1)^n\,\hpg{2}{1}{-n,\;-a-n}{\!1-a+m-n}{\frac1{1-z}}\\
&=&\frac{m!\,(a+1)_n}{(n+m)!}\,z^n\,\hpg{2}{1}{-n,\;m+1}{a+1}{\,1-\frac1{z}\,}\\
&=&\frac{m!\,(a+1)_n}{(n+m)!}\;\hpg{2}{1}{-n,\;a-m}{a+1}{\,1-z\,}.
\end{eqnarray*}
For the last four expressions, standard Euler's formula (\ref{flinear1}) can
be applied; this gives us 4 non-terminating hypergeometric expressions. In
total we have 6 terminating and 4 non-terminating hypergeometric expressions
for each basis solution in (\ref{abbasis}). The remaining 4 Kummer's
solutions are related by Euler-Pfaff transformations; they represent one
Gauss hypergeometric function. The relation between this non-terminating and
the two terminating solutions is a consequence of
formula (\ref{althpgrel}) in Lemma \ref{f21ff}:%\vspace{-2.3pt}
\begin{eqnarray} \label{abnonterm}
&&\hspace{-30pt}(1-z)^{-a}\;\hpg{2}{1}{-m,-a-n}{-n-m}{\,z}=\hpg{2}{1}{-n,a-m}{-n-m}{\,z}
+\nonumber\\
&&\quad(-1)^m\,\frac{n!\,m!\,(a-m)_{n+m+1}}{(n+m)!\,(n+m+1)!}\,z^{n+m+1}\,
\hpg{2}{1}{m+1,\,a+n+1}{n+m+2}{\,z}.
\end{eqnarray}

\section{The trivial monodromy group}
\label{trivial}

Here we consider hypergeometric equations with the trivial monodromy group.
That means that solutions can be meromorphically continued to the entire
projective line $\PP^1$, so they  are rational functions.
By part 5 of Theorem \ref{dencases}, general hypergeometric equation with
trivial monodromy group is $E(-n,\ell+1,-m-n)$. We have the following three
terminating solutions:%\vspace{-2.3pt}
\begin{eqnarray} \label{threesols}
\hpg21{-n,\,\ell+1}{-n-m}{\,z},\qquad
(1-z)^{-\ell-1}\,\hpg21{-m,\,\ell+1}{-n-m}{\frac{z}{z-1}},\nonumber\\
z^{n+m+1}\,(1-z)^{-m-\ell-1}\,\hpg21{-\ell,\,n+1}{-m-\ell}{1-z}.
\end{eqnarray}
Each of them can be transformed by formulas
(\ref{termseries1})--(\ref{termseries6}). For example,%\vspace{-2.2pt}
\begin{eqnarray*}
\hpg21{-n,\,\ell+1}{-n-m}z &=& (1-z)^n\;\hpg{2}{1}{-n,-n-m-\ell-1}{-n-m}{\frac{z}{z-1}}\\
&=&\frac{m!\,(n+\ell)!}{\ell!\,(n+m)!}\,z^n\,\hpg{2}{1}{-n,\;m+1}{-n-\ell}{\,\frac1{z}\,}\\
&=&\frac{m!\,(n+\ell)!}{\ell!\,(n+m)!}\,(z-1)^n\,\hpg{2}{1}{-n,-n-m-\ell-1}{-n-\ell}{\frac1{1-z}}\\
&=&\frac{m!\,(n+m+\ell+1)!}{(n+m)!\,(m+\ell+1)!}\,z^n\,\hpg{2}{1}{-n,\;m+1}{2+m+\ell}{\,1-\frac1{z}\,}\\
&=&\frac{m!\,(n+m+\ell+1)!}{(n+m)!\,(m+\ell+1)!}\;\hpg{2}{1}{-n,\;\ell+1}{2+m+\ell}{\,1-z\,}.
\end{eqnarray*}
Note how permutation of the numbers $m,n,\ell$ permutes the three sets of
hypergeometric representations for solutions in (\ref{threesols}), if we
ignore the front factors and change of the argument. The last two identities
can be transformed to non-terminating series by Euler's formula, so in total
we have 6 terminating and 2 non-terminating hypergeometric expressions for
each of the three solutions. This exhausts the 24 Kummer's series. Relation
between the three solutions is a consequence of (\ref{abnonterm}):
\begin{eqnarray*}
&&\hspace{-30pt}(1-z)^{-\ell-1}\;\hpg{2}{1}{-m,\ell+1}{-n-m}{\frac{z}{z-1}}
=\hpg{2}{1}{-n,\ell+1}{-n-m}{\,z\,}+\nonumber\\
&&\quad(-1)^m\,\frac{n!\,(m+\ell)!}{\ell!\,(n+m)!}\,z^{n+m+1}\,(1-z)^{-m-\ell-1}\,
\hpg{2}{1}{-\ell,\,n+1}{-m-\ell}{1-z}.
\end{eqnarray*}

\section{Additive monodromy group}
\label{additive}

Here we consider hypergeometric equations whose monodromy group is (up to
conjugation) a non-trivial subgroup of $\GG_a$. By part 6 of Theorem
\ref{dencases}, general hypergeometric equation of this type is
$E(-\ell,-n-\ell,-m-n-2\ell)$. The point $z=0$ is not a logarithmic point
for this equation, there is a basis of power series solutions there:
\begin{equation} \label{additivbasis}
\hpg21{-\ell,-n-\ell}{-m-n-2\ell}z,\qquad
z^{m+n+2\ell+1}\,\hpg21{m+\ell+1,m+n+\ell+1}{m+n+2\ell+2}z.
\end{equation}
The first solution has the following terminating expressions:
\begin{eqnarray*}
\hpg21{-\ell,-n-\ell}{-n-m-2\ell}z
&=&(1-z)^{\ell}\;\hpg{2}{1}{-\ell,-m-\ell}{-n-m-2\ell}{\frac{z}{z-1}}\\
&=&C_1\,(-z)^\ell\,\hpg{2}{1}{-\ell,\;n+m+\ell+1}{n+1}{\,\frac1{z}\,}\\
&=&C_1\,(1-z)^\ell\,\hpg{2}{1}{-\ell,-m-\ell}{n+1}{\frac1{1-z}}\\
&=&C_2\,z^\ell\,\hpg{2}{1}{-\ell,\,n+m+\ell+1}{m+1}{\,1-\frac1{z}\,}\\
&=&C_2\;\hpg{2}{1}{-\ell,\;-n-\ell}{m+1}{\,1-z\,}\\
&=& (1-z)^{-m}\;\hpg{2}{1}{-m-\ell,-n-m-\ell}{-n-m-2\ell}{\,z\,}\\
&=&(1-z)^{n+\ell}\;\hpg{2}{1}{-n-\ell,-n-m-\ell}{-n-m-2\ell}{\frac{z}{z-1}}\\
&=&C_1\,(-z)^{m+\ell}(1-z)^{-m}\hpg{2}{1}{-m-\ell,\;n+\ell+1}{n+1}{\,\frac1{z}\,}\\
&=&C_2\,z^{n+\ell}\,\hpg{2}{1}{-n-\ell,\,m+\ell+1}{m+1}{\,1-\frac1{z}\,},\\
\end{eqnarray*}
where
\[
C_1=\frac{(n+\ell)!\,(n+m+\ell)!}{n!\,(n+m+2\ell)!},\qquad
C_2=\frac{(m+\ell)!\,(n+m+\ell)!}{m!\,(n+m+2\ell)!}.
\]
This solution has also non-terminating hypergeometric expressions, with the
argument $1-z$ or $1/(1-z)$ %$\frac1{1-z}$
by Euler's formula. In total we have 10 terminating and 2 non-terminating
hypergeometric expressions for this solution. The number of distinct
terminating expressions may drop to 8 (if $m=0$ or $n=0$) or to 6 (if
$m=n=0$). The second solution in (\ref{additivbasis}) has 2, 3 or 4 distinct
hypergeometric expressions due to Euler-Pfaff transformations.

Other Kummer's series at $z=1$ and $z=\infty$ are undefined (or coincide
with terminating expressions, if $m=0$ or $n=0$). Consequently, there is no
basis of power series solutions at these points $z=1$, $z=\infty$; they are
logarithmic. % for the equation $E(-\ell,-n-\ell,-m-n-2\ell)$.
In the following Theorem, we present logarithmic expressions for the
function
\begin{equation}\label{u3exp0}
U_3=\frac{(-1)^{m+1}}{(m+n+2\ell+1)!}\,z^{m+n+2\ell+1}\,
\hpg{2}{1}{m\!+\!\ell\!+\!1,\,m\!+\!n\!+\!\ell\!+\!1}{\!m+n+2\ell+2}{\,z\,}.
\end{equation}
Notice that all terms with the $\psi$-function can be written as terminating
sums of rational numbers, and that all sums in expression (\ref{u3exp1}) are
terminating.
\begin{theorem} \label{loglogterm}
Set $C_3=1\big/\ell!\,(m+\ell)!\,(n+\ell)!\,(m+n+\ell)!$. The following
formulas hold:
%\[C_3=\frac1{\ell!\,(m+\ell)!\,(n+\ell)!\,(m+n+\ell)!}.\]
\begin{eqnarray} \label{u3exp1}
%U_3&\equal&\frac{(-1)^{m+1}}{(m+n+2\ell+1)!}\,z^{m+n+2\ell+1}\,
%\hpg{2}{1}{m\!+\!\ell\!+\!1,\,m\!+\!n\!+\!\ell\!+\!1}{\!m+n+2\ell+2}{\,z\,}\\
U_3&\equal&C_3\,(m+n+2\ell)!\;\hpg21{-\ell,\,-n-\ell}{-m-n-2\ell}{\,z\,}\log(1-z)\nonumber\\
&&+\sum_{k=0}^\ell\frac{\psi(n\!+\!\ell\!-\!k\!+\!1)+\psi(\ell\!-\!k\!+\!1)
-\psi(m\!+\!k\!+\!1)-\psi(k\!+\!1)}{(m+k)!\,(n+\ell-k)!\,(\ell-k)!\,k!}(1-z)^k\nonumber\\
&&-(z-1)^{-m}\,%+(-1)^{m+1}(1-z)^{-m}
\sum_{k=0}^{m-1}\frac{(m-k-1)!}{(m+n+\ell-k)!\,(m+\ell-k)!k!}\,(z-1)^k\nonumber\\
&&+(-1)^{\ell}(z-1)^{n+\ell} %+(-1)^n\,(1-z)^{n+\ell}\,
\sum_{k=0}^{n-1}\frac{(n-k-1)!}{(m+n+\ell-k)!\,(n+\ell-k)!\,k!}\,\frac1{(z-1)^k}\\
%U_3&\equal&\frac{(-1)^{m+1}}{(m+n+2\ell+1)!}\,z^{m+n+2\ell+1}\,
%\hpg{2}{1}{m\!+\!\ell\!+\!1,\,m\!+\!n\!+\!\ell\!+\!1}{\!m+n+2\ell+2}{\,z\,}\\
%&\equal&C_3\,(m+n+2\ell)!\;\hpg21{-\ell,\,-n-\ell}{-m-n-2\ell}{\,z\,}\log(1-z)\nonumber\\
%&&+\sum_{k=0}^\ell\frac{\psi(n\!+\!\ell\!-\!k\!+\!1)+\psi(\ell\!-\!k\!+\!1)
%-\psi(m\!+\!k\!+\!1)-\psi(k\!+\!1)}{(m+k)!\,(n+\ell-k)!\,(\ell-k)!\,k!}(1\!-\!z)^k\nonumber\\
%&&-(z-1)^{-m}\,\sum_{k=0}^{m-1}\frac{(m-k-1)!}{(m+n+\ell-k)!\,(m+\ell-k)!k!}\,(z-1)^k\nonumber\\
%&&-(z\!-\!1)^{\ell+1}\,\sum_{k=0}^{n-1}\!\frac{k!}{(n\!-\!k\!-\!1)!(m\!+\!\ell\!+\!k\!+\!1)!
%(\ell\!+\!k\!+\!1)!}\,(z-1)^k\\
\label{u3exp2}&\equal&C_3\,(m+n+2\ell)!\,\hpg21{-\ell,\,-n-\ell}{-m-n-2\ell}{\,z\,}\log(1-z)\nonumber\\
&&-C_3\,z^{m+n+2\ell+1}\,(z-1)^{-m}\,\sum_{k=0}^{m-1}\frac{(n\!+\!\ell\!+\!k)!\,(\ell\!+\!k)!\,
(m\!-\!k\!-\!1)!}{k!}\,(z-1)^k\nonumber\\
&&+C_3\,z^{m+n+2\ell+1}\,\sum_{k=0}^\infty\frac{(m+\ell+k)!\,(m+n+\ell+k)!}
{(m+k)!\,k!}\,(1-z)^k\,\times\nonumber\\
&&\hspace{20pt}\big(\psi(m\!+\!n\!+\!\ell\!+\!k\!+\!1)+\psi(m\!+\!\ell\!+\!k\!+\!1)
-\psi(m\!+\!k\!+\!1)-\psi(k\!+\!1)\big)\\
%&\equal&C_3\,(m+n+2\ell)!\,\hpg21{-\ell,\,-n-\ell}{-m-n-2\ell}{\,z\,}\log\frac1{1-z}\nonumber\\
%&&+(-1)^{n+1}C_3\,z^{m+n+2\ell+1}\,(1-z)^{n+\ell}\,\sum_{k=0}^{n-1}\frac{(m\!+\!\ell\!+\!k)!\,(\ell\!+\!k)!\,
%(n\!-\!k\!-\!1)!}{k!}\,\frac1{(z-1)^k}\nonumber\\
%&&+C_3\,z^{m+n+2\ell+1}\,(1-z)^{\ell}\sum_{k=0}^\infty\frac{(n+\ell+k)!\,(m+n+\ell+k)!}
%{(n+k)!\,k!}\,\frac1{(1-z)^k}\,\times\nonumber\\
%&&\hspace{20pt}\big(\psi(m\!+\!n\!+\!\ell\!+\!k\!+\!1)+\psi(n\!+\!\ell\!+\!k\!+\!1)
%-\psi(n\!+\!k\!+\!1)-\psi(k\!+\!1)\big)\\
\label{u3exp3}&\equal&C_3\,(m+n+2\ell)!\;\hpg21{-\ell,\,-n-\ell}{-m-n-2\ell}{\,z\,}\,\log\frac{1-z}{z}\nonumber\\
&&-\frac{z^{m+\ell}\,(z-1)^{-m}}{(n+\ell)!(m+n+\ell)!}\,\sum_{k=0}^{m-1}
\frac{(n+\ell+k)!\,(m-k-1)!}{(m+\ell-k)!\,k!}\,\frac{(z-1)^k}{z^k}\nonumber\\
&&\hspace{-1pt}+\frac{z^\ell}{(n+\ell)!(m+n+\ell)!}\,\sum_{k=0}^\ell\!\frac{(m+n+\ell+k)!}
{(m+k)!(\ell-k)!\,k!} \times\nonumber\vspace{-1pt}\\ &&\hspace{24pt}
\left(\psi(m\!+\!n\!+\!\ell\!+\!k\!+\!1)+\psi(\ell\!-\!k\!+\!1)-\psi(m\!+\!k\!+\!1)-\psi(k\!+\!1)\right)
\frac{(z-1)^k}{z^k}\nonumber\\
&&+\frac{(-z)^\ell}{(n+\ell)!(m+n+\ell)!}\,\sum_{k=\ell+1}^{\infty}\,
\frac{(m\!+\!n\!+\!\ell\!+\!k)!(k-\ell-1)!}{(m+k)!\,k!}\,\frac{(z-1)^k}{z^k}\\
\label{u3exp4}&\equal&C_3\,(m+n+2\ell)!\;\hpg21{-\ell,\,-n-\ell}{-m-n-2\ell}{\,z\,}\,\log\frac{1-z}{z}\nonumber\\
&&-\frac{z^{m+n+\ell}\,(z-1)^{-m}}{\ell!\,(m+\ell)!}\,\sum_{k=0}^{m-1}
\frac{(\ell+k)!\,(m-k-1)!}{(m+n+\ell-k)!\,k!}\,\frac{(z-1)^k}{z^k}\nonumber\\
&&\hspace{-1pt}+\frac{z^{n+\ell}}{\ell!\,(m+\ell)!}\,\sum_{k=0}^{n+\ell}\!\frac{(m+\ell+k)!}
{(m+k)!(n+\ell-k)!\,k!} \times\nonumber\vspace{-1pt}\\ &&\hspace{24pt}
\left(\psi(m\!+\!\ell\!+\!k\!+\!1)+\psi(n\!+\!\ell\!-\!k\!+\!1)-\psi(m\!+\!k\!+\!1)-\psi(k\!+\!1)\right)
\frac{(1-z)^k}{z^k}\nonumber\\
&&+\frac{(-z)^{n+\ell}}{\ell!\,(m+\ell)!}\,\sum_{k=n+\ell+1}^{\infty}\,
\frac{(m\!+\!\ell\!+\!k)!(k-n-\ell-1)!}{(m+k)!\,k!}\,\frac{(z-1)^k}{z^k}.
\end{eqnarray}
Besides, in each expression one can interchange $m$ and $n$, provided that
$z$ is replaced by $z/(z-1)$ and the whole expression is multiplied by
$(-1)(1-z)^\ell$.
\end{theorem}
\proof To derive the formulas, we consider the equation
$E(-\ell,-n-\ell,-m-n-2\ell)$ in the context of Theorem \ref{logtermt}. Then
we have $U_3=U_2\big/m!(n+\ell)!\,\ell!$.
In formula (\ref{u2exp3}), we get rid of the singular $\psi$-values and the
tangent term by using Lemma \ref{psitransf}. The result is:
\begin{eqnarray*}
U_3&\equal&\frac1{m!\,(n+\ell)!\,\ell!}\,\hpg21{-\ell,\,-n-\ell}{m+1}{1-z}\,\log(1-z)\nonumber\\
&&+(-1)^{m+1}\,\,(1-z)^{-m}\,\sum_{k=0}^{m-1}
\frac{(m-k-1)!}{(m+n+\ell-k)!\,(m+\ell-k)!k!}\,(z-1)^k\nonumber\\
&&+\sum_{k=0}^\ell\frac{\psi(n\!+\!\ell\!-\!k\!+\!1)+\psi(\ell\!-\!k\!+\!1)-\psi(m\!+\!k\!+\!1)-\psi(k\!+\!1)}
{(m+k)!\,(n+\ell-k)!\,(\ell-k)!\,k!}(z\!-\!1)^k\nonumber\\
&&+(-1)^\ell\,\sum_{k=\ell+1}^{\ell+n}\!\frac{(k-\ell-1)!}{(n+\ell-k)!(m+k)!\,k!}\,(z-1)^k.
\end{eqnarray*}
We apply Euler's transformation (\ref{flinear1}) to the first hypergeometric
sum, rewrite the last sum in the opposite direction, and get (\ref{u3exp1}).
Formulas (\ref{u3exp2}) and (\ref{u3exp3}) are just rewritten expressions
(\ref{u2exp4}) and (\ref{u2exp6}), respectively. Formula (\ref{u3exp4}) can
be obtained from (\ref{u2exp6}) after interchanging the first two parameters
of $E(-\ell,-n-\ell,-m-n-2\ell)$; the same formula can be obtained by
carefully applying Lemma \ref{psitransf} to expression (\ref{u2exp5}).

To see the last statement, one can check the described transformation on
formula (\ref{u3exp0}) and compare it with Pfaff's transformation. (Formula
(\ref{u3exp1}) is invariant under this transformation as well, up to Euler's
transformation of the first term and summing the second term in the opposite
direction. Interchanging the singular points $z=1$, $z=\infty$ produces the
same transformation.)\qed

\bibliographystyle{alpha}
\bibliography{../hypergeometric}

\end{document}